\documentclass [a4paper,english]{article}
\usepackage{amsmath, amsfonts, amssymb}
\usepackage{babel,indentfirst,epigraph}
\usepackage[dvips]{graphicx,color}

\title{The Damascus Inequality}
\author{Fozi M. Dannan, Sergey M. Sitnik}

\date{}

\begin{document}
\maketitle

\section{The problem formulation}

In 2016 Prof. Fozi M. Dannan from Damascus, Syria proposed the next inequality

\begin{equation}\label{1}
\frac{x-1}{x^2-x+1}+\frac{y-1}{y^2-y+1}+\frac{z-1}{z^2-z+1} \leq 0,
\end{equation}
providing that $xyz=1$ for $x,y,z > 0.$
It became widely known but was not proved yet in spite of elementary formulation.

An obvious generalization is the next inequality

\begin{equation}\label{N}
\sum_{k=1}^{n} \frac{x_k-1}{x_k^2-x_k+1} \le 0,
\end{equation}
providing that $x_1\cdot x_2 \ldots \cdot x_n=1$ for $x_k \ge 0, 1\le k \le n$.

It is obvious that (\ref{N}) is true for $n=1$, it is easy to prove it also for $n=2$ directly.
But it is not true for $n=4$ as follows from an example with $x_1=x_2=x_3=2, x_4=\frac{1}{8}$, then
(\ref{N}) is reducing to $1-\frac{56}{57}\le 0$ which is untrue.

As a consequence   (\ref{N}) is also untrue for any $n\ge 4$ due to an example with $x_1=x_2=x_3=2, x_4=\frac{1}{8}, x_5=\ldots =x_n=1$. So the only non--trivial case in (\ref{N}) is $n=3$.

In this paper we prove inequality \eqref{1} together with similiar ones
\begin{eqnarray}
\frac{1}{x^2-x+1}+\ \frac{1}{y^2-y+1}+\frac{1}{z^2-z+1}\ \ \le 3 \label{2}\\
\frac{x}{x^2-x+1}+\ \frac{y}{y^2-y+1}+\frac{z}{z^2-z+1}\ \ \le 3 \label{3}\\
\frac{x-1}{x^2+x+1}+\ \frac{y-1}{y^2+y+1}+\frac{z-1}{z^2+z+1}\ \ \le 0 \label{4}\\
\frac{1}{x^2+x+1}+\ \frac{1}{y^2+y+1}+\frac{1}{z^2+z+1}\ \ \ge 1 \label{5}\\
\frac{x}{x^2+x+1}+\ \frac{y}{y^2+y+1}+\frac{z}{z^2+z+1}\ \ \le 1 \label{6}\\
\frac{x+1}{x^2+x+1}+\ \frac{y+1}{y^2+y+1}+\frac{z+1}{z^2+z+1}\ \ \le 2. \label{7}
\end{eqnarray}

Also some generalizations will be considered.

\section{Proof of the main inequality \eqref{1}}

\textbf{Theorem 1.} An inequality  \eqref{1} holds true providing that $xyz=1$ for $x,y,z > 0.$

For the proof we  need an auxiliary inequality that seems to be very interesting by itself.

\textbf{ Lemma 1.}  Let $x, y, z\ $ be positive numbers such that  $xyz=1$.  Then
\begin{equation}\label{lem}
x^2+y^2+z^2-3\left(x+y+z\right)+6\ge 0
\end{equation}
holds true.

Note that inequality \eqref{lem} is not a consequence of well--known family of Klamkin--type inequalities for symmetric functions \cite{MPF}. So \eqref{lem} is a new quadratic Klamkin--type inequality in three variables under restriction $xyz=1$.
Due to its importance we give three  proofs to it based on  different ideas.

\textbf{First proof of Lemma 1.}

To prove (\ref{lem}) let introduce the Lagrange function
$$
L(x,y,z,\lambda)= x^2 + y^2 + z^2  - 3 (x + y + z) +6 - \lambda(xyz-1).
$$
On differentiating it follows
$$
\lambda=x^2-2x=y^2-2y=z^2-2z.
$$
It follows that at the minimum (it obviously exists) $x=y$, so three variables at the minimum are $x,y=x,z=1/x^2$. From $x^2-2x=z^2-2z$ we derive the equation in $x$:
$$
2x^2-3x=2/x^4-3/x^2, f(x)=2x^6-3x^5+3x^2-2=0.
$$
One root is obvious $x=1$. Let us prove that there are no other roots for $x\ge 0$. Check that derivative is positive
$$
f'(x)=12x^5-15x^4+6x=3x(4x^4-5x^3+2)\ge 0, x\ge 0.
$$
Define a function $g(x)=4x^4-5x^3+2$, its derivative $g'(x)$ has one zero for $x\ge 0$ at $x=15/16$ and the function $g(x)$ is positive at this zero at its minimum $g(15/16)=15893/16384>0$. So $g(x)$ is positive, $f(x)$ is strictly increasing on $x\ge 0$, so $f(1)=0$ is its only zero.

\textbf{Second proof of Lemma 1.}

Consider the function
\[f\left(x,y\right)=x^2+y^2+\frac{1}{x^2y^2}-3\left(x+y+\frac{1}{xy}\right)+6\ \ ,\]
where $x, y, z $ are positive numbers.
 We show that $f(x,y)$ attains its minimum  0 at $x=1, y=1$ using partial derivative test.

Calculate
\begin{equation} \label{10}
\frac{\partial f}{\partial x}=2x-\frac{2}{x^3y^2}-3+\frac{3}{x^2y}=0,
\end{equation}
\begin{equation}  \label{11}
\frac{\partial f}{\partial y}=2y-\frac{2}{x^2y^3}-3+\frac{3}{xy^2}=0.
\end{equation}
Now multiplying  \eqref{10} and \eqref{11} respectively by  $x$ and $-y$ and adding to obtain
\[\left(x-y\right)\left(2x+2y-3\right)=0\ \ .\]

Here we have two cases.

Case I.  $x=y\ $, which implies from equation \eqref{10} that
\[2x^6-2-3x^5+3x^2=0\]
or
\begin{equation} \label{12}
\left(x^3-1\right)\left(2x^3-3x^2+2\right)=0.
\end{equation}
Equation \eqref{12} has only one positive root $x=1$ and consequently $y=1.$ Notice that the equation
\[2x^3-3x^2+2=0\]
does not have  positive roots because for $x\ge 0$
 the function
   $u\left(x\right)=2x^3-3x^2+2$
    satisfies the following properties :
    $\left(i\right) u\left(0\right)=2, \left(ii\right)\ \mathrm{min} u\left(x\right)=1$
     at $x=1$,
    $\left(iii\right) u\left(\infty \right)= \infty.$
Therefore $f\left(x,y\right)$ attains its maximum or minimum at $x=1,\ y=1.$

 Case II. $2x+2y=3\ \ .$  Adding \eqref{10} and \eqref{11} we get

\[2\left(x+y\right)-6-2(\frac{1}{x^3y^2}+\frac{1}{x^2y^3})+\ 3\left(\frac{1}{x^2y}+\frac{1}{xy^2}\right)=0\]
or
\[-3-\frac{3}{x^3y^3}+\frac{3}{2x^2y^2}=0\]
and
\[-6x^3y^3+3xy-6=0\ \ .\]
Putting $t=xy$  we obtain
\[2t^3-t+2=0 .\]
In fact this equation does not have positive root (notice that $t=xy$   should be positive). This is because the function $u=2t^3-t+2\ $ satisfies the following properties :
\[\left(i\right)\ \ u\left(0\right)=2,\ \ \]
\[\left(ii\right)\ \mathrm{for}\ \ t>0, {\mathrm{min}\,u\left(t\right)=u\left(\frac{1}{\sqrt{6}}\right)>0},\ \left(iii\right)\ \ \ u\left(\infty \right)>0.\]
The last step is to show that
\[f\left(x,y\right)\ge f\left(1,1\right)=0.\]
It is enough to show that $f\left(s,t\right)>f\left(1,1\right)$  for at least one point $(s,t)\neq (1,1)$  .  Take for example $f\left(2,3\right)=\frac{7}{2}+\frac{1}{36}$.

\textbf{Third proof of lemma 1 (Geometrical Method).}

Geometrically it is enough to prove that the surface $xyz=1$ lies outside the sphere $(x-3/2)^2+(y-3/2)^2+(z-3/2)^2=3/4$ except the only intersection point (1,1,1) as it is shown on the next graph:

\begin{center}
\includegraphics[scale=0.8]{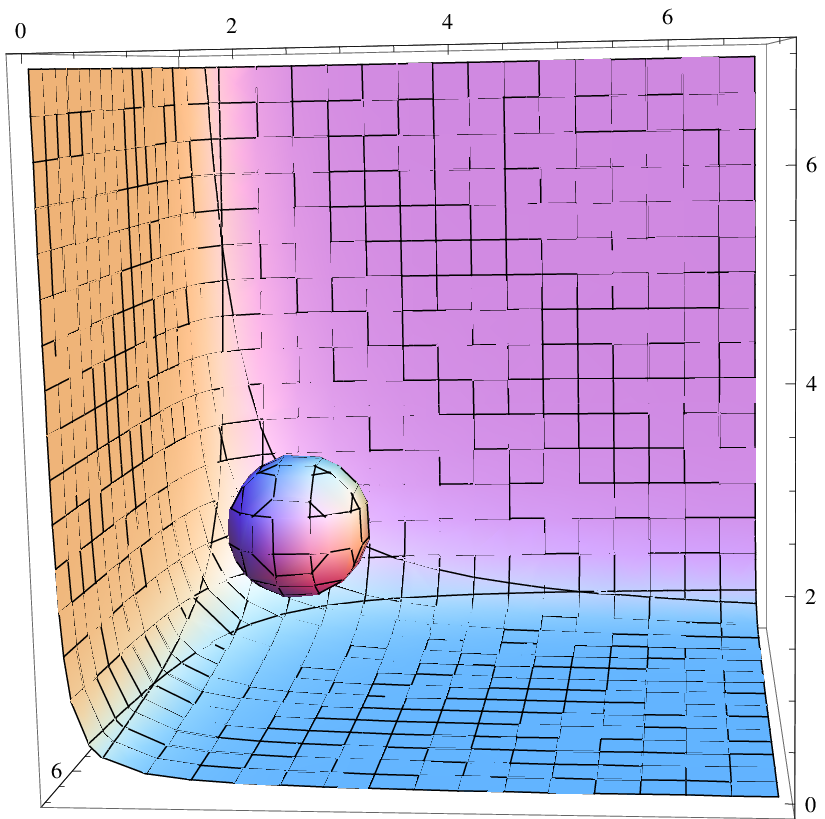}
\end{center}

Let $M$ and $S$  be surfaces defined by

$M:\ \ xyz=1$  and $S:\ {(x-\frac{3}{2})}^2+{(y-\frac{3}{2})}^2+{(z-\frac{3}{2})}^2-\frac{3}{4}\ =0$

 1. If   ${(z-\frac{3}{2})}^2-\frac{3}{4}\ge 0\ \ $ then
\[{(x-\frac{3}{2})}^2+{(y-\frac{3}{2})}^2+{(z-\frac{3}{2})}^2-\frac{3}{4}\ge 0\ \ \]
and equivalently
\[x^2+y^2+z^2-3\left(x+y+z\right)+6\ge 0.\ \ \ \]

2. If  ${(z-\frac{3}{2})}^2-\frac{3}{4}\le 0\ \ $ then $\frac{3\ -\ \sqrt{3}}{2}\le z\le \frac{3\ +\ \sqrt{3}}{2}.$

3. We take  horizontal sections for both  $M$  and so get for any plane
\[\frac{3\ -\ \sqrt{3}}{2}\le z=k\le \frac{3\ +\ \sqrt{3}}{2}\]
two curves: equilateral hyperbola $H(k)$ with vertex $(\frac{1}{\sqrt{k} }, \frac{1}{\sqrt{k}})$   and a circle  $C\left(k\right)$  which radius is given by
\[r^2\left(k\right)={\frac{3}{4}-(k-\frac{3}{2})}^2=-k^2+3k-\frac{3}{2}\ \ \]
with center  at $(\frac{3}{2}, \frac{3}{2}, k)$.

4.  For $z=1\ $, we have the hyperbola  $xy=1$  and the circle
\[{(x-\frac{3}{2})}^2+{(y-\frac{3}{2})}^2=\frac{1}{2}.\]

5.  We show that the distance $d(v,c)$  between the vertex of the hyperbola and the center  $\ c\ $ of the circle is always greater than or equal to the radius of the circle. The distance  $d(v,c)$   is given by
\[{d^2\left(v,c\right)=\left(\frac{1}{\sqrt{k}\ }-\frac{3}{2}\right)}^2+{\left(\frac{1}{\sqrt{k}\ }-\frac{3}{2}\right)}^2\mathrm{=2\ }{\left(\frac{1}{\sqrt{k}\ }-\frac{3}{2}\right)}^2.\]
The radius is given by
\[r^2\left(k\right)={\frac{3}{4}-(k-\frac{3}{2})}^2=-k^2+3k-\frac{3}{2}\ \ \ .\]
We need to show that the vertex is always outside the circle i.e.
 $d^2\left(v,c\right)\ge r^2\left(k\right)$  for all
\[\frac{3\ -\ \sqrt{3}}{2}\le k\le \frac{3\ +\ \sqrt{3}}{2}\ \ .\ \ \]
Clearly that  $d\left(v,c\right)=r$  for $k=1$  and the hyperbola tangents the circle at the point
$\left(1, 1, 1\right)$ .

 For  $\frac{3\ -\ \sqrt{3}}{2}\le k<1$  , as  $k\ $ decreases from  $1\ $ to   $\frac{3\ -\ \sqrt{3}}{2}$ , the radius of the circle becomes smaller . From the other side the vertex  $(\frac{1}{\sqrt{k}\ }, \frac{1}{\sqrt{k}}, k)$ moves away from  $(1, 1, 1)$  towards a point $ \left(0, 0, k \right)$. This follows from the distance function of the vertex
\[Ov=\ \frac{\sqrt{2}}{\sqrt{k}}, \left(\ 0<k_1\le k_2<1\to \ \frac{\sqrt{2}}{\sqrt{k_2}}<\frac{\sqrt{2}}{\sqrt{k_1}}\ \ \right).\]

6. For  $1<k\le \frac{3\ +\ \sqrt{3}}{2}$   , we show that
\[{d^2\left(v,c\right)=g\left(k\right)=2\ \left(\frac{1}{\sqrt{k}\ }-\frac{3}{2}\right)}^2>r^2={\frac{3}{4}-(k-\frac{3}{2})}^2=h\left(k\right).\]
In fact, $h\left(k\right)$ is a concave down parabola and has its maximum at  $k=1\ $, i.e.
$\mathrm{max}\  h\left(k\right)=\frac{1}{2}$ and $h\left(k\right)$ is decreasing for $k>1.$  Also,  $g\left(1\right)=\frac{1}{2}$  and $g(k)$ is increasing for  $k>1$ because
$g'(k)=4 \left(-\frac{1}{k\sqrt{k}}\right)\left(\frac{1}{\sqrt{k}}-\frac{3}{2}\right)>0$  for
$k \ge 1$ (notice that  $\frac{1}{\sqrt{k}}<1$).
Therefore  $g\left(k\right)>h\left(k\right)$   for  $k>1$  . Eventually we conclude that
\[{(x-\frac{3}{2})}^2+{(y-\frac{3}{2})}^2+{(z-\frac{3}{2})}^2-\frac{3}{4}\ge 0\ \ \]
and consequently
\[x^2+y^2+z^2-3\left(x+y+z\right)+6\ge 0\ \ \ \ \ \]
 for all values of  $\left(x, y, z \right)$  that satisfy  $xyz=1$.

\textbf{Proof of the theorem 1.}

Now consider the inequality to prove \eqref{1}.
After simplifying with the use of Wolfram Mathematica it reduces to
\begin{eqnarray*}
-3 + 3 x - 2 x^2 + 3 y - 3 x y + 2 x^2 y - 2 y^2 + 2 x y^2 -
 x^2 y^2 + 3 z - 3 x z + 2 x^2 z-\\ - 3 y z + 3 x y z - 2 x^2 y z +
 +2 y^2 z - 2 x y^2 z + x^2 y^2 z - 2 z^2 + 2 x z^2 - x^2 z^2 +\\+
 2 y z^2 - 2 x y z^2 + x^2 y z^2 - y^2 z^2 + x y^2 z^2\le 0.
\end{eqnarray*}
Using \textbf{SymmetricReduction} function of Wolfram Mathematica we derive
\begin{eqnarray*}
3 - x y - x z - y z + 3 x y z - 3 (x + y + z) + 2 (x + y + z)^2 -
 x y z (x y + x z + y z)\\ -
 -2 (x + y + z) (x y + x z + y z) + (x y + x z + y z)^2
\ge 0.
\end{eqnarray*}
Using $xyz=1$ let further simplify
\begin{eqnarray*}
6 - 3 (x + y + z) + 2 (x + y + z)^2 - 2 (x y + x z + y z) -\\
 - 2 (x + y + z) (x y + x z + y z) + (x y + x z + y z)^2\ge 0.
\end{eqnarray*}
In terms of elementary symmetric functions
$$
S_1=x+y+z,  S_2=xy+yz+xz
$$
it is
\begin{equation}\label{K1}
S_2^2 -2S_1 S_2 -2 S_2 +2S_1^2 -3 S_1 +6 \ge 0.
\end{equation}
As $S_2^2 -2S_1 S_2+S_1^2\ge 0$ it is enough to prove
\begin{equation}\label{K2}
S_1^2 -2 S_2  -3 S_1 +6 \ge 0.
\end{equation}

Expanding it again in $x,y,z$ we derive an inequality to prove for positive variables
 \begin{equation}\label{Eq}
x^2 + y^2 + z^2  - 3 (x + y + z) +6 \ge 0.
\end{equation}
But this is exactly an inequality from Lemma 1. So Theorem 1 is proved.

\section{Proof of  inequalities \eqref{2}--\eqref{7}}

Let us start with   two propositions.

\textbf{Proposition 1.}

For any real numbers $u,\ v,\ w$ such that
\[\left(1+u\right)\left(1+v\right)\left(1+w\right)>0,\]
the inequality
\[\frac{1}{1+u}+\frac{1}{1+v}+\frac{1}{1+w}\le k\ \ \ \ \ \ (\ge k)\]
is equivalent to
\[kuvw+\left(k-1\right)\left(uv+vw+wu\right)+\left(k-2\right)\left(u+v+w\right)+
k-3 \ge 0\ \  (\le 0). \]

\textbf{Proposition 2.}

For any real numbers $u\ ,\ v\ ,w$ such that
\[\left(u-1\right)\left(v-1\right)\left(w-1\right)>0\]
the inequality
\[\frac{1}{u-1}+\frac{1}{v-1}+\frac{1}{w-1}\le k\ \ \ \ \ \ (\ge k)\]
is equivalent to
\[kuvw-\left(k+1\right)\left(uv+vw+wu\right)+\left(k+2\right)\left(u+v+w\right)
-(k+3) \ge 0 \ \ (\le 0). \]

The validity of propositions 1 and 2 can be obtained by direct expansions.

\textbf{Proof that $\eqref{1}\Leftrightarrow \eqref{2}$.}

In fact
\[\frac{x-1}{x^2-x+1}+\ \frac{y-1}{y^2-y+1}+\frac{z-1}{z^2-z+1}=\]
\[=\frac{x^2-(x^2-x+1)}{x^2-x+1}+\ \frac{y^2-(y^2-y+1)}{y^2-y+1}+\frac{z^2-(z^2-z+1)}{z^2-z+1}=\]
\[=-3+\sum_{cyc}{\frac{x^2}{x^2-x+1}}.\]
Now if the right side is $\le 0$  then
\[\sum_{cyc}{\frac{x^2}{x^2-x+1}}\le 3\]
and  consequently
\[\sum_{cyc}{\frac{1}{1-\left(1/x\right)+{(1/x)}^2}}\le 3.\]

\textbf{Proof of 5.}

We need to prove
\[\sum_{cyc}{\frac{1}{x^2+x+1}}\ge 1\]
Let $u=x^2+x\ \ ,\ v=y^2+y\ ,\ \ w=z^2+z$  .  Using Proposition 1 the required inequality can be written as follows :
\[uvw-\left(u+v+w\right)-2\le 0 .\]
Going back to $x,y,z$  we get
\[\left(x+1\right)\left(y+1\right)\left(z+1\right)-\left(x^2+y^2+z^2\right)-\left(x+y+z\right)-2\le 0 .\]
Or
\[xy+yz+zx\le x^2+y^2+z^2\]
which is obvious.

\textbf{Proof of 6.}

It follows from elementary calculus that for any real number $x$  we have
\[\frac{x}{x^2+x+1}\le \frac{1}{3}\]
and the inequality follows directly.

\textbf{Proof that $\eqref{5}+\eqref{6} \Rightarrow \eqref{4}$.}

Really adding together \eqref{6} with \eqref{5} multiplied by $-1$ we derive \eqref{4}.

\textbf{Proof that  $\eqref{5} \Rightarrow \eqref{7}$.}

The required inequality is equivalent to
\[\sum_{cyc}{\frac{x^2+x+1-x^2}{x^2+x+1}}=3-\sum_{cyc}{\frac{x^2}{x^2+x+1}\ }\le 2\]
or
\[\sum_{cyc}{\frac{x^2}{x^2+x+1}\ }\ge 1\]
which is true from inequality \eqref{5}.

\section{Modifications of original inequality}

In this section we consider modifications of the original inequality (\ref{1}) providing that $xyz=1$ for $x,y,z \ge 0.$

1. An inequality (\ref{1}) is equivalent to

\begin{equation}\label{M1}
\frac{x^2-1}{x^3-1}+\frac{y^2-1}{y^3-1}+\frac{z^2-1}{z^3-1} \leq 0.
\end{equation}

This form leads to generalization with more powers, cf. below.

2. An inequality (\ref{1}) is equivalent to

\begin{equation}\label{M2}
\frac{x^2}{x^2-x+1}+\frac{y^2}{y^2-y+1}+\frac{z^2}{z^2-z+1} \leq 3.
\end{equation}

3. Let take $x\rightarrow \frac{1}{x}, y\rightarrow \frac{1}{y}, z\rightarrow \frac{1}{z}$. Then we derive another equivalent form of the inequality (\ref{1})

\begin{equation}\label{M3}
\frac{x-x^2}{x^2-x+1}+\frac{y-y^2}{y^2-y+1}+\frac{z-z^2}{z^2-z+1} \leq 0,
\end{equation}
due to the functional equation
\begin{equation}\label{FE}
f(\frac{1}{x})=-x f(x)
\end{equation}
for the function
\begin{equation}\label{fun}
f(x)= \frac{x-1}{x^2-x+1}.
\end{equation}

So it seems possible to generalize the original inequality in terms of functional equations too.

To one more similar variant leads a change of variables $x\rightarrow xy, y\rightarrow yz, z\rightarrow xz$:
\begin{equation}\label{M4}
\frac{xy-1}{x^2 y^2-xy+1}+\frac{yz-1}{y^2 z^2-yz+1}+\frac{xz-1}{x^2z^2-xz+1} \leq 0,
\end{equation}
or like (\ref{M3})
\begin{equation}\label{M5}
\frac{xy-x^2y^2}{x^2 y^2-xy+1}+\frac{yz-y^2z^2}{y^2 z^2-yz+1}+\frac{xz-x^2z^2}{x^2z^2-xz+1} \leq 0,
\end{equation}

It is also possible to consider generalizations of (\ref{1}) under the most general transformations $x\rightarrow g(x,y,z), y\rightarrow h(x,y,z), z\rightarrow \frac{1}{g(x,y,z)h(x,y,z)}$ with positive functions
$g(x,y,z), h(x,y,z)$ still preserving a condition $xyz=1$.

4. A number of cyclic inequalities follow from previous ones by a substitution
$$
x=\frac{a}{b},  y=\frac{b}{c},   z=\frac{c}{a}, xyz=1.
$$
On this way we derive from \eqref{1}, \eqref{2}--\eqref{7} the next cyclic inequalities:
\begin{eqnarray}
\frac{ab-b^2}{a^2-ab+b^2}+\frac{bc-c^2}{b^2-bc+c^2}+\frac{ca-a^2}{c^2-ca+a^2} \leq 0,\\
\frac{b^2}{a^2-ab+b^2}+\frac{c^2}{b^2-bc+c^2}+\frac{a^2}{c^2-ca+a^2}\ \ \le 3 \\
\frac{ab}{a^2-ab+b^2}+\frac{bc}{b^2-bc+c^2}+\frac{ca}{c^2-ca+a^2}\ \ \le 3 \\
\frac{ab-b^2}{a^2+ab+b^2}+\frac{bc-c^2}{b^2+bc+c^2}+\frac{ca-a^2}{c^2+ca+a^2}\ \ \le 0 \\
\frac{b^2}{a^2+ab+b^2}+\frac{c^2}{b^2+bc+c^2}+\frac{a^2}{c^2+ca+a^2}\ \ \ge 1 \\
\frac{ab}{a^2+ab+b^2}+\frac{bc}{b^2+bc+c^2}+\frac{ca}{c^2+ca+a^2}\ \ \le 1 \\
\frac{ab+b^2}{a^2+ab+b^2}+\frac{bc+c^2}{b^2+bc+c^2}+\frac{ca+a^2}{c^2+ca+a^2}\ \ \le 2.
\end{eqnarray}

On cyclic inequalities among which  Schur, Nessbit and Shapiro ones are the most well--known cf. \cite{MPF}--\cite{MOa}.

5. Some geometrical quantities connected with trigonometric functions and triangle geometry satisfy a condition $xyz=1$, cf. \cite{MPV}--\cite{M}. For example, we may use in standard notations for triangular geometry values:
$$
x=\frac{a}{4p}, \ \ y=\frac{b}{R},\ \ z=\frac{c}{r}  ;
$$
$$
x=\frac{a+b}{2}, \ \ y=\frac{b+c}{p},\ \ z=\frac{a+c}{p^2+r^2+2rR}  ;
$$
$$
x=R h_a, \ \ y=\frac{h_b}{2p^2},\ \ z=\frac{h_c}{r^2}  ;
$$
$$
x=2R^2 \sin(\alpha), \ \ y=\frac{\sin(\beta)}{r},\ \ z=\frac{\sin(\gamma)}{p}  ;
$$
$$
x=(p^2-4R^2-4rR-r^2)\tan(\alpha), \ \ y=\frac{\tan(\beta)}{2p},\ \ z=\frac{\tan(\gamma)}{r}  ;
$$
$$
x=\frac{\tan(\alpha)}{\tan(\alpha)+\tan(\beta)+\tan(\gamma)}, y=\frac{\tan(\beta)}{\tan(\alpha)+\tan(\beta)+\tan(\gamma)}, z=\frac{\tan(\gamma)}{\tan(\alpha)+\tan(\beta)+\tan(\gamma)};
$$
$$
x=\tan(\alpha/2), \ \ y=p \tan(\beta/2),\ \ z=\frac{\tan(\gamma/2)}{r}  ;
$$
$$
x=\frac{a}{4(p-a)}, \ \ y=\frac{b}{R(p-b)},\ \ z=r \frac{c}{p-c}  ;
$$
$$
x=4R \sin(\alpha/2), \ \ y=\sin(\beta/2),\ \ z=\frac{\sin(\gamma/2)}{r}  ;
$$
$$
x=4R \cos(\alpha/2), \ \ y=\cos(\beta/2),\ \ z=\frac{\cos(\gamma/2)}{p} .
$$

6. The  above geometrical identities of the type $x\,y\,z=1$ which we use for applications of considered inequalities are mostly consequences of Vieta's formulas \cite{SM}. It is interesting to use these formulas for cubic equation directly.

\textbf{Theorem 2.} Let $x,y,z$ be positive roots of the cubic equation with any real $a,b$
$$
t^3 + a t^2 + b t -1=0.
$$
The  for these roots $x,y,z$ all inequalities of this paper are valid.

7. We can generalize inequalities \eqref{2}, \eqref{5}--\eqref{7} for more general powers.
For this aim we use Bernoulli's inequalities \cite{MPF}--\cite{Mv} :  for  $u>0\ $
the following inequalities hold true
\[u^{\alpha }-\alpha u+\alpha -1\ge 0\ ,\ \ (\ \alpha >1\ \ \mathrm{or}\ \alpha <0),\]
\[u^{\alpha }-\alpha u+\alpha -1\le 0\ ,\ \ \left(\ 0<\alpha <1\ \ \right)\ .\]

\textbf{Lemma 2. }
Assume that $x, y, z $ are positive numbers such that $xyz=1.$
Then the following inequality holds true :

\[{\left(\frac{1}{x^2-x+1}\right)}^{\alpha }+{\left(\frac{1}{y^2-y+1}\right)}^{\alpha }+{\left(\frac{1}{z^2-z+1}\right)}^{\alpha }\ \ \le 3\ \ \ \ \ \]
for   $0<\alpha <1.$

\textbf{Proof}. Let
\[X=x^2-x+1\ ,\ Y=y^2-y+1\ ,\ Z=\ z^2-z+1\ \ .\ \ \]
Then we have
\[{\left(\frac{1}{x^2-x+1}\right)}^{\alpha }+{\left(\frac{1}{y^2-y+1}\right)}^{\alpha }+{\left(\frac{1}{z^2-z+1}\right)}^{\alpha }\ \ \le \ \ \]
\[\le \alpha \left(\frac{1}{X}+\frac{1}{Y}+\frac{1}{Z}\right)+3\left(1-\alpha \right)\le 3\ \ \ .\ \ \ \ \]

Similarly we have from \eqref{6} that
\[{\left(\frac{x}{x^2+x+1}\right)}^{\alpha }+{\left(\frac{y}{y^2+y+1}\right)}^{\alpha }+{\left(\frac{z}{z^2+z+1}\right)}^{\alpha }\le 3-2\alpha \ \ \ \]
and from \eqref{7} we have
\[{\left(\frac{x+1}{x^2+x+1}\right)}^{\alpha }+{\left(\frac{y+1}{y^2+y+1}\right)}^{\alpha }+{\left(\frac{z+1}{z^2+z+1}\right)}^{\alpha }\le 3-\alpha \ .\ \ \ \]
For   $\alpha >1$  or $\alpha <0$ we have from \eqref{5}

\[{\left(\frac{1}{x^2+x+1}\right)}^{\alpha }+{\left(\frac{1}{y^2+y+1}\right)}^{\alpha }+{\left(\frac{1}{z^2+z+1}\right)}^{\alpha }\ge 3-2\alpha.\ \ \]

\section{Generalizations of original inequality to ones with a set of restrictions on symmetric functions}

It is easy to show that the maximum of the function (\ref{fun}) is attained for $x\ge 0$ at $x=2$ and equals to 1/3.
\begin{center}
\includegraphics[scale=0.4]{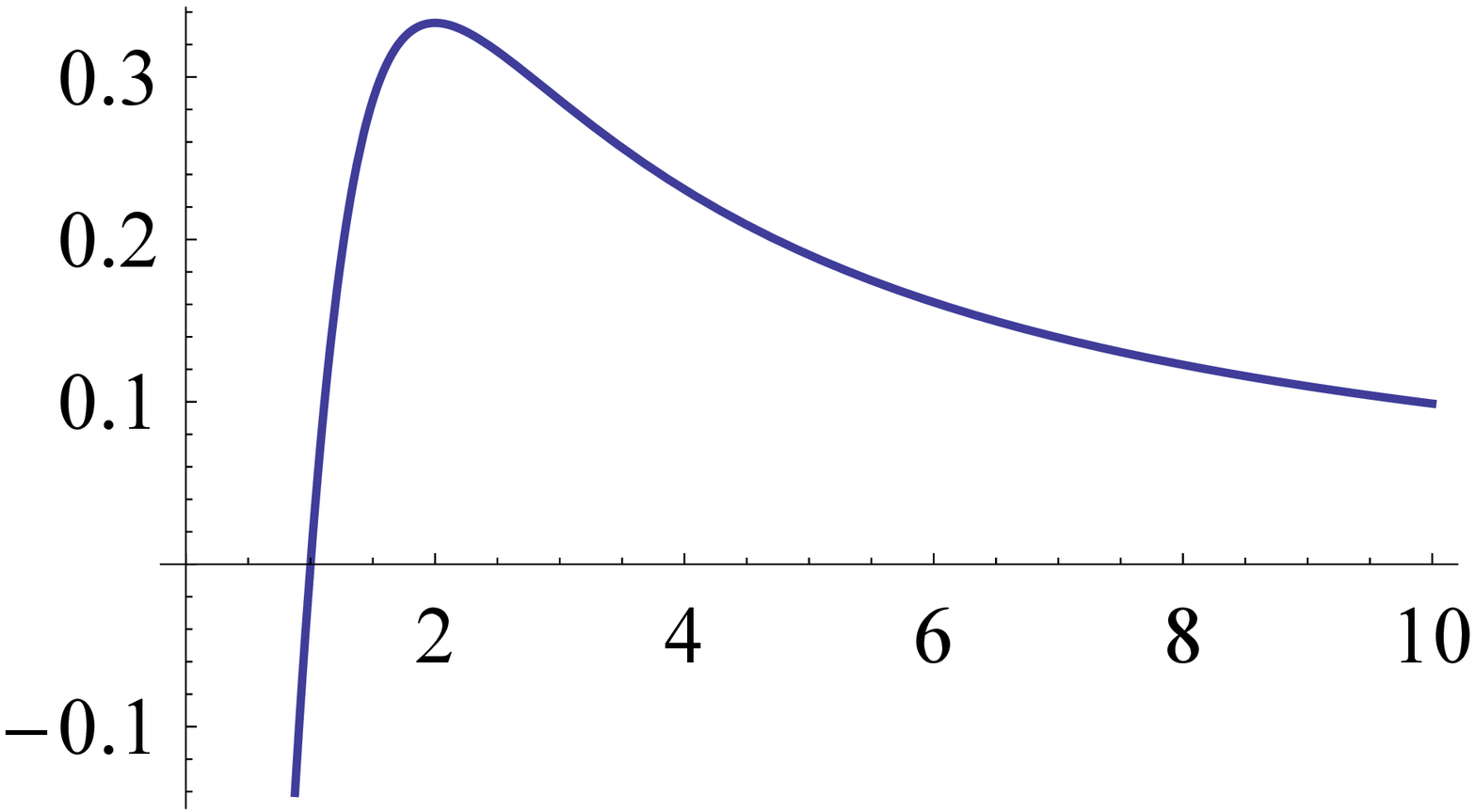}
\end{center}

So the next unconditional inequality holds
\begin{equation}\label{N2}
\sum_{k=1}^{k=n} \frac{x_k-1}{x_k^2-x_k+1} \le \frac{n}{3}; \ \ \ x_k \ge0
\end{equation}

Consider symmetric functions
$$
S_1=\sum_{k=1}^{k=n} x_k, S_2=\sum_{\substack{k,m=1,\\ k\neq m}}^{n} x_k \cdot x_m,\dots, S_n=x_1x_2\cdots x_n.
$$

\textbf{The generalized Damascus inequality}

Prove an inequality
\begin{equation}\label{GDI}
\sum_{k=1}^{k=n} \frac{x_k-1}{x_k^2-x_k+1} \le \frac{n}{3}-C(a_1,a_2,\cdots,a_n ); \ \ \ x_k \ge0
\end{equation}
and find the best positive constant in it under conditions on symmetric functions
\begin{equation}\label{res}
S_1=a_1, S_2=a_2, \cdots, S_n=a_n
\end{equation}
with may be some restrictions in (\ref{res}) omitted.

The unconditional constant for positive numbers in (\ref{GDI}) is $C=0$ and the original inequality gives
$C=\frac{n}{3}$ in case $n=3$ and a single restriction $S_3=1$ in the list (\ref{res}).

It seems that a problem to find the sharp constant in the inequality (\ref{GDI}) under general conditions (\ref{res}) is a difficult problem.

For three numbers so more inequalities of the type (\ref{GDI}) may be considered, e.g.

1. Prove inequality \eqref{GDI} for positive numbers under condition $S_1=1$ and find the best constant for this case.

2. Prove inequality \eqref{GDI} for positive numbers under condition $S_2=1$ and find the best constant for this case.

Also combined conditions may be considered.

3. Prove inequality \eqref{GDI} for positive numbers under conditions $S_1=a, S_2=b$ and find the best constant $C(a,b)$ in \eqref{GDI} for this case.

\section{Symmetricity of  symmetric inequalities}

 There are many inequalities that are written in terms of symmetric functions  as  $F\left(p,q\right)\ \le 0\ \ \ \ \left(\ \ge 0\ \ \right)$  , where
\[p=S_1=x+y+z,\ q=S_2=xy+yz+zx,\ \  r=S_3=xyz=1.\]
The following Lemma enlarge the amount of inequalities that one can obtain as a series of very complicated inequalities.

\textbf{Lemma 3.}  If the inequality
\[F\left(p,q\right)\ \le 0\ \ \ \ \left(\ \ge 0\ \ \right)\]
holds true , then the following inequalities are satisfied :
\[\left(i\right)\ \ \ \ \ \ F\left(q,p\right)\ \le 0\ \ \ \ (\ \ \ge 0\ \ \ ),\]
and
\[\left(ii\right)\ \ \ F\left(\ \ q^2-2p\ ,\ \ p^2-2q\right)\ \le 0\ \ \ \ \ \ \left(\ \ge 0\ \ \right)\ .\]

\textbf{Proof.} (i). Assume that
\[F(p,q)\ = F(x+y+z,xy+yz+zx)\ge   0.\]
Using  transformations
\[x\to xy, y\to yz, z\to zx\]
we obtain
\[F\left(p,q\right)=F\left(xy+yz+zx,\ xyyz+yzzx+zxxy\right)=\]
\[=F\left(xy+yz+zx,\ x+y+z\ \right)=F\left(q,p\right)\ge 0.\]

Notice that we can also use  transformations
\[x\to \frac{1}{x}, y\to \frac{1}{y}, z\to \frac{1}{z}.\]

(ii). Now assume that
\[F(p,q)\ = F(x+y+z,xy+yz+zx)\ge   0.\]
 Using  transformations
\[x\to \frac{xy\ }{z},  y\to \frac{yz}{x}, z\to \frac{zx}{y}\]
we derive
\[\frac{xy\ }{z}+\ \frac{yz}{x}\ \ \ +\frac{zx}{y}\ =\ x^2y^2+y^2z^2+z^2x^2=\ \]
\[={(xy+yz+zx)}^2-2\left(x+y+z\right)=q^2-2p\ .\ \ \]
Also it follows
\[\frac{xy\ }{z}\ \frac{yz}{x}\ \ \ +\ \frac{yz}{x}\ \frac{zx}{y}\ +\ \frac{zx}{y}\ \frac{xy\ }{z}=\ \]
\[=\ \frac{y\ }{zx}+\ \frac{z}{xy}\ \ \ +\frac{x}{yz}=\ x^2+y^2+z^2\ \ \]
\[={\left(x+y+z\right)}^2-2\left(xy+yz+zx\ \right)=p^2-2q\ \ .\]

The proof is complete.

At the end we propose an unsolved problem.

\textbf{Problem.} Find all possible non--negative values of four variables $x_1,x_2,x_3,x_4$ with restriction $x_1\cdot x_2\cdot x_3\cdot x_4=1$ for which the next inequality holds
\begin{equation}\label{N4}
\sum_{k=1}^{4} \frac{x_k-1}{x_k^2-x_k+1} \le 0,
\end{equation}
As we know from the example at the beginning of the paper the inequality \eqref{N4} is not true for all such values, e.g. it fails for $x_1=x_2=x_3=2,x_4=1/8$.

\vskip 10mm

AUTHORS:\\
\\
Fozi M. Dannan,\\
Department of Basic Sciences, Arab International University,
P.O.Box  10409, Damascus , SYRIA,
e-mail : fmdan@scs-net.org
\\
\\
Sergei M. Sitnik,\\
Voronezh Institute of the Ministry of Internal Affairs of Russia,
Voronezh, Russia,
e-mail : pochtaname@gmail.com
\end{document}